\documentclass{aptpub}

\authornames{Didier Piau}
\shorttitle{Asymptotics of posteriors for binary branching processes}

\begin{document}

\title{Asymptotics of posteriors for binary branching processes}

\authorone[Universit\'e Joseph Fourier]{Didier Piau}

\addressone{Institut Fourier UMR 5582, 
Universit\'e Joseph Fourier Grenoble 1,
100 rue des Maths, BP 74, 38402 Saint Martin d'H\`eres, France.
E-mail address: \texttt{Didier.Piau@ujf-grenoble.fr}.
Webpage: \texttt{www-fourier.ujf-grenoble.fr/\%7edpiau}.
}

\begin{abstract}
  We compute the posterior distributions of the initial
  population and parameter of binary branching processes, in the limit
  of a large number of generations. We compare this Bayesian procedure
  with a more na\"\i ve one, based on hitting times of some random
  walks. In both cases, central limit theorems are available, with
  explicit variances.
\end{abstract}



%
\section{Introduction}
\label{s.intro}
This paper is devoted to some estimation procedures of binary branching
processes in a Bayesian setting.  To be more specific, 
let $(X_n)_{n\ge0}$ denote a Galton-Watson
process which starts from the initial population
$X_0\ge1$ and whose offspring is ruled by the distribution
$$
(1-U)\,\delta_1+U\delta_2\quad\mbox{with}\ 0<U<1,
$$  
where $\delta_x$ denotes the Dirac mass at $x$. This means that, at
every generation, each individual dies and is replaced by $1$ or $2$
individuals, with probability $1-U$ and $U$ respectively,
independently of the fate of the other individuals, and that $X_n$
counts generation $n$.

In a Bayesian framework, the initial population $X_0$
and the offspring parameter $U$ are both random and unknown. To keep things simple, we
also assume that $X_0$ and $U$ are independent, and we wish to estimate them
from the observation of a finite path $x_{1:n}=(x_k)_{1\le k\le n}$ 
of the process $X_{1:n}=(X_k)_{1\le k\le n}$
up to
a given time $n\ge1$.

Well known motivations for such a study are various biological
settings where one observes $X_{1:n}$ but $X_0$ and $U$ are unknown.
One example is the modeling of polymerase chain reaction.
Probabilistic models of polymerase chain reactions were proposed and
studied by Sun (1995), Weiss and von Haeseler (1995) and (1997),
Peccoud and Jacob (1996), Piau (2002), (2004), (2005), and Jagers and
Klebaner (2003). Recently, Lalam and Jacob (2007) introduced and
studied the Bayesian setting above, see also Lalam (2007).  For other
Bayesian approaches of branching processes, see Scott (1987), Prakasa
Rao (1992), Mendoza and Guti\'errez-Pe\~na (2000), and, for the
interesting model of bisexual branching process, Molina, Gonz\'alez
and Mota (1998) for example. Finally, the idea of studying a branching
process backwards, but to estimate its age rather than its initial
population, is in Klebaner and Sagitov~(2002).

In models of polymerase chain
reactions and in similar contexts, the
initial population $X_0$ is the size of a small sample,
extracted at random from a much larger population.
This suggests that the initial population 
$X_0$ should be Poisson distributed, say with
parameter $\Lambda$. We assume that $\Lambda$ is random as well.
Jeffreys' principle, see Kass and Wasserman~(1996), 
then indicates that the prior distributions of
$\Lambda$ and $U$ should be proportional to measures which we compute
below.  To sum up the result of these computations,
the prior of $\Lambda$ is easy to write down but
improper and the prior of $U$ is awkward but proper.  However, the
posterior of $(X_0,U)$ conditionally on $X_{1:n}$
is a proper distribution, which can be computed explicitly.  
In particular, this posterior distribution depends only on 
$X_1$, $X_n$ and $S_n=X_1+\cdots+X_n$. Unfortunately, 
it is also rather unwieldy.

In such situations, one may rely on numerical algorithms, 
based on MCMC for example, to
simulate the posterior distributions with any prescribed degree of
accuracy. Rather, we look for simple asymptotics in realistic
regimes.  Namely, we assume that $n$ is large and we are interested in
the asymptotic posterior distribution of $(X_0,U)$ assuming that $X_n$
is large and that the ratio $S_n/X_n$ converges to a finite limit.
This assumption is almost surely fulfilled by the paths of
binary branching processes since these are supercritical. In
this setting, we show that the posterior distributions indeed converge
and we compute explicitly their limit.

\section{Results}

To describe our results, we introduce some notations. 
Let  $x_{0:\infty}=(x_n)_{n\ge0}$ denote a sequence of positive 
integers. We say that such a sequence is admissible if, 
for every nonnegative $n$,
$x_n\le x_{n+1}\le2x_n.$
We say that an admissible sequence is regular if furthermore, $x_n/s_n$
converges to a positive limit when $n$ goes to infinity, 
where $s_n=x_1+\cdots+x_n$.
The binary index $B(x_{0:\infty})$ of a 
regular admissible sequence $x_{0:\infty}$ is the real number in $]0,1]$
defined by 
$$
B(x_{0:\infty})=\lim_{n\to\infty}\frac{x_{n+1}}{s_n}.
$$
The renormalized index $R(x_{0:\infty})$ of a 
regular admissible sequence $x_{0:\infty}$  is the real number 
in $[0,+\infty[$ defined by
$$
R(x_{0:\infty})=\lim_{n\to\infty}\frac{(s_n-x_{n+1})^2}{4x_{n+1}s_n}.
$$
Almost every (sequence which can be
realized as a) path of a binary branching process is admissible
and regular.
The renormalized index is a function of the 
binary index, namely
$R(x_{0:\infty})=\varrho(B(x_{0:\infty}))$ 
where, for every $u$ in $]0,1]$,
$$
\varrho(u)=\frac{(1-u)^2}{4u}.
$$
The binary index and the normalized index are asymptotic quantities, in
the sense that, for every
nonnegative integer $n$, the indexes of a regular admissible sequence
$x_{0:\infty}$ do not depend on the first values $x_{0:n}$.

From now on, letters $k$ and $n$ are used to enumerate generations of the
process (that is, the time) and symbols $x$, $x_k$, $x_n$ and $y$ 
are used to measure population sizes.

\begin{definition}[Distributions]
\label{d.munu}
  For every positive real number $r$ and every positive
  integer $x$, the finite discrete measure $\nu(r,x)$ and the
  discrete probability measure $\mu(r,x)$, both on the positive
  integers, are defined by
$$
\nu(r,x)=\sum_{y=h(x)}^x\binom{2y}{y}\,\binom{y}{x-y}\,r^y\,\delta_y,
\qquad
\mu(r,x)=\frac{\nu(r,x)}{|\nu(r,x)|}.
$$
For every positive integer $x$, the integer $h(x)$ in the formula above is 
the upper half of $x$, that is, the 
smallest integer  such that $2h(x)\ge x$. In other words,
$h(2x)=h(2x-1)=x$ for every positive integer $x$. 
\end{definition}

Our main result is as follows. 

\begin{theorem}[Posterior distributions]
\label{ta}
(1)
The path $X_{0:\infty}$ of a binary branching process with parameter $U$ 
is almost
surely regular admissible and its binary index is almost surely 
 $B(X_{0:\infty})=U$.  

(2)
Assume that the prior distribution of $(X_0,U)$ satisfies Jeffreys' principle. 
Then, for every regular admissible sequence $x_{1:\infty}$ with binary index
$u=B(x_{1:\infty})$ in $]0,1[$,
 the posterior distribution of
$(X_0,U)$ conditionally on $X_{1:n}=x_{1:n}$
converges when $n$ goes to infinity
to the distribution
$
\mu(\varrho(u),x_1)\otimes\delta_{u}.
$
\end{theorem}

Theorem~\ref{ta} shows that the limit posterior distribution of $X_0$
when $n$ goes to infinity is almost surely $\mu(r,x)$ with
$r=\varrho(U)$ and $x=X_1$. Unless $r=0$, $r=1$ or $x=1$, $\mu(r,x)$
is not degenerate, hence the value of $X_0$ can be determined only
with some uncertainty, even from an infinite trajectory
$X_{1:\infty}$. On the contrary, $U$ is a function of the 
infinite trajectory $X_{1:\infty}$.

The limit distribution $\mu(\varrho(u),x_1)\otimes\delta_{u}$ in
theorem~\ref{ta} converges to the Dirac distribution at $(x_1,0)$ when
$u$ converges to $0$ and to the Dirac distribution at $(h(x_1),1)$
when $u$ converges to $1$.  Our next result describes the intuitively
obvious variations of $\mu(r,x)$ with respect to $r$ and $x$. First,
since $r=\varrho(u)$ is a decreasing function of $u$ and the offspring
distribution of the branching process is stochastically increasing
with $u$, one should expect $\mu(r,x)$ to increase stochastically when
$r$ increases. Likewise, since $x$ represents the population at time
$1$, one should expect $\mu(r,x)$, which represents the population at
time $0$, to increase stochastically when $x$ increases.

We recall that a measure $\mu_1$ is stochastically larger than a
measure $\mu_2$ if and only if 
$\mu_1([z,+\infty))\ge\mu_2([z,+\infty))$ for every 
real number $z$. 

\begin{proposition}[Ordering of limit posterior distributions]
\label{lo}
For every positive integer $x$, the family $(\mu(r,x))_{r\ge0}$ is
stochastically increasing. For every positive real number $r$,
the family $(\mu(r,x))_{x\ge1}$ is stochastically increasing.
\end{proposition}

We now characterize the limit of
$\mu(r,x)$ for every fixed value of $r$, when $x$ converges to infinity.
 
\begin{theorem}[Limit posterior distributions of initial populations]
\label{tb}
  Fix $u$ in $]0,1[$. For every positive integer $x$, let $\xi_x$ 
denote a random variable with distribution $\mu(\varrho(u),x)$. When
  $x$ converges to infinity, 
the expectation and the mode of $\xi_x/x$
  both converge to 
$$m_u=1/(1+u),$$
and
the random variables
  $\left(\xi_x-m_ux\right)/\sqrt{x}$ converge in distribution
  to a centered Gaussian distribution with variance 
$$
\sigma^2_u=u(1-u)/(1+u)^3.
$$
\end{theorem}

For the sake of comparison, we turn to another natural way to estimate
initial populations of branching processes with known offspring
distributions, based on hitting times.  To describe this in the
setting of binary branching processes, we first introduce some notations.

\begin{definition}[Hitting times]
Fix a real number $u$ in $]0,1[$, and let $(\varepsilon_x)_{x\ge1}$ denote
a sequence of independent Bernoulli random variables with distribution
$(1-u)\,\delta_1+u\,\delta_2$. For every positive integer $x$, let
$\sigma_x:=\varepsilon_1+\cdots+\varepsilon_x$.
Define the distribution of the hitting time 
$\eta_x$ by the relation
$$
\mathbb{P}(\eta_x=y)=\mathbb{P}(\sigma_y=x\,|\,H_x),\qquad
\mbox{where}\ H_x=\{\exists z\ge1\,;\,\sigma_z=x\}.
$$
\end{definition}

When the value of $u$ is known, an 
estimation procedure of $X_0$ based on $X_1=x$ is to propose
the value $y$ for $X_0$ with probability $\mathbb{P}(\eta_x=y)$, 
thus an estimator of $X_0$ when $X_1=x$ is the distribution of $\eta_x$.

Recall that $m_u=1/(1+u)$ and $\sigma^2_u=u(1-u)/(1+u)^3$.

\begin{theorem}[Initial populations through hitting times]
\label{tc}
Fix a real number $u$ in $]0,1[$. For every positive integer $x$,
$$
\left|\mathbb{E}(\eta_x)-m_ux\right|\le2u/(1+u)^2\le1/2.
$$
Furthermore, when $x$ converges to infinity,
$(\eta_x-m_ux)/\sqrt{x}$ converges in distribution to a centered
Gaussian variable with variance 
$\sigma^2_u$.  
\end{theorem}

The rest of the paper is organized as follows. We prove
theorem~\ref{ta} and proposition~\ref{lo} in section~\ref{s.ta}
and theorem~\ref{tb} in section~\ref{s.tb}.
Finally, the proof of theorem~\ref{tc}, sharper bounds on $\mathbb{E}(\eta_x)$
and a brief comparison with another, non Bayesian, estimation procedure 
are in section~\ref{s.tc}.

%
\section{Posterior distributions}
\label{s.ta}
\subsection{Preliminaries}

Jeffreys' principle, see Kass and Wasserman~(1996), indicates that the
prior measure for a parameter $\theta$ governing the distribution
$\nu_\theta$ of a random variable $Z$ should have a density
proportional to $J(\theta)^{1/2}$, where
$$
J(\theta)=-\mathbb{E}_{\theta}\left(\frac{\partial^2}{\partial\theta^2}\log\nu_\theta(Z)\right).
$$
We apply this to the parameter $(\Lambda,U)$. Parts of lemma~\ref{l.prior} are in Lalam and Jacob~(2007).

\begin{lemma}\label{l.prior}
  For every positive integer $n$, the prior measure for $(\Lambda,U)$ according to Jeffreys' principle and
based on $X_{0:n}$
  is the product of the prior measures for $\Lambda$ and $U$. 
The prior
  measures for $\Lambda$ and for $U$ are respectively proportional to
the measures  
$\mathrm{d}\lambda/\sqrt\lambda$ on $\lambda>0$ and $\pi_n(u)\,\mathrm{d} u$ on  $0<u<1$,
  where
$$
\pi_n(u)=\sqrt{\frac{(1+u)^n-1}{u^2(1-u)}}.
$$
In particular, the prior of $U$ is proper.
\end{lemma}

\begin{proof}[Proof of lemma \ref{l.prior}] Assume that $X_0$ is Poisson
  distributed with parameter $\Lambda$ and that $X_{0:n}$ is a binary
  branching process with parameter $U$. Then the distribution
  $\nu_{\Lambda,U}$ of $X_{0:n}$ is such that
$$
\nu_{\Lambda,U}(x_{0:n})=\mathrm{e}^{-\Lambda}\frac{\Lambda^{x_0}}{x_0!}\prod_{k=1}^n\binom{x_{k-1}}{x_k-x_{k-1}}
U^{x_k-x_{k-1}}(1-U)^{2x_{k-1}-x_k}.
$$
Up to a factor $C(x_{0:n})$ which does not depend on $(\Lambda,U)$,
$\log\nu_{\Lambda,U}(x_{0:n})$
is
$$
-\Lambda+x_0\log\Lambda+
(x_n-x_0)\log U+(s_n-2x_n+2x_0)\log(1-U)+C(x_{0:n}).
$$
This is the sum of a function of $\Lambda$ and a function of $U$, hence the 
prior measures are product measures.
As regards the prior for $\Lambda$, 
$$
\frac{\partial^2}{\partial\Lambda^2}\log\nu_\Lambda(x_0)
=
-\frac{x_0}{\Lambda^2},
\qquad\mbox{hence}\
J(\Lambda)=\frac{\mathbb{E}_\Lambda(X_0)}{\Lambda^2}=\frac1{\Lambda}.
$$
As regards the prior for $U$,
$$\frac{\partial^2}{\partial U^2}\log\nu_U(x_{0:n})
=
-(x_n-x_0)/U^2-(s_n-2x_n+2x_0)/(1-U)^2,
$$
hence
$$
J_n(U)=\frac{\mathbb{E}_U(X_n-X_0)}{U^2}+\frac{\mathbb{E}_U(S_n-2X_n+2X_0)}{(1-U)^2},
$$
where $S_n=X_1+\cdots+X_n$.  Since
$\mathbb{E}_U(X_k)=(1+U)^k\mathbb{E}(X_0)$ for every nonnegative
integer $k$, one finds that $J_n(U)=\mathbb{E}(X_0)\pi_n(U)^2$ with
the notations of the lemma.

Finally, up to multiplicative constants, $\pi_n(u)$ behaves like
$1/\sqrt{u}$ when $u$ converges to $0$ and like $1/\sqrt{1-u}$ when
$u$ converges to $1$. Hence, $\pi_n$ is integrable and there exists a
(proper) prior distribution for $U$.  This concludes the proof of lemma \ref{l.prior}.
\end{proof}

From now on, we fix a positive integer $n$, 
we assume that the observations are
$X_{1:n}=x_{1:n}$ with $x_{1:n}=(x_k)_{1\le k\le n}$ and we recall
that $s_n=x_1+\cdots+x_n$. The posterior distribution in
lemma~\ref{l5} is similar, but not equal, to a posterior distribution
computed in Lalam and Jacob~(2007).

\begin{lemma}\label{l5}
  The posterior distribution of $(X_0,U)$ conditionally on
  $X_{1:n}=x_{1:n}$ depends only on $x_1$, $x_n$ and $s_n$, and
  is proportional to the measure
$$
\sum_{x=h(x_1)}^{x_1}2^{-2x}\,\binom{2x}{x}\,\binom{x}{x_1-x}\,
u^{x_n-x}\,(1-u)^{s_n-2x_n+2x}\,\pi_n(u)\,\delta_x\otimes\mathrm{d} u.
$$
\end{lemma}

\begin{proof}[Proof of lemma \ref{l5}]
  Fix $u$, $x_{1:n}$ and $x$ such that $h(x_1)\le x\le x_1$. Then, the
  conditional probability $\mathbb{P}(U\in\mathrm{d}
  u,X_0=x\,\vert\,X_{1:n}=x_{1:n})$ is proportional to
$$
\nu_U(\mathrm{d} u)\int\nu_\Lambda(\mathrm{d}\lambda)
\mathbb{P}_\lambda(X_0=x)\mathbb{P}_u(X_{1:n}=x_{1:n}\,\vert\,X_0=x),
$$
where $\nu_U(\mathrm{d} u)=\pi_n(u)\mathrm{d} u$ and  
$\nu_\Lambda(\mathrm{d}\lambda)=\mathrm{d}\lambda/\sqrt\lambda$.
Hence,
$$
\int\nu_\Lambda(\mathrm{d}\lambda)\mathbb{P}_\lambda(X_0=x)
=
\frac{\Gamma(x+1/2)}{\Gamma(x+1)}
=
\sqrt\pi\,2^{-2x}\,\binom{2x}{x}.
$$
Likewise, using the computations in the proof of lemma \ref{l.prior}, one gets
$$
\mathbb{P}_u(X_{1:n}=x_{1:n}\,\vert\,X_0=x)=C(x_{1:n})\,\binom{x}{x_1-x}\,
u^{x_n-x}\,(1-u)^{s_n-2x_n+2x},
$$
where $C(x_{1:n})$ does not depend on $(x,u)$. This concludes the proof of lemma \ref{l5}.
\end{proof}

\subsection{Proof of theorem \ref{ta}}
Part (1) follows from the fact that, when $n$ converges to infinity,
 $X_n/(1+U)^n$ converges almost surely to a random positive and finite limit.

A sketch of the proof of part (2) is as follows. 
Consider the distribution in lemma~\ref{l5} 
and assume that $x_n$ converges to infinity and
that $x_n/(s_n-x_n)$ converges to $v$. Then $s_n-x_n$ is equivalent to 
$x_n/v$, hence 
$$
u^{x_n}(1-u)^{s_n-2x_n}=\left(u^v(1-u)^{1-v}\right)^{s_n-x_n+o(x_n)}.
$$
The inner parenthesis is maximal when $u=v$, and the exponent
converges to infinity, hence this contribution becomes concentrated
around the value $u=v$.  The remaining factor involving $u$ 
in the distribution described in lemma \ref{l5}
is $\varrho(u)^x$, and the convergence to $\mu(\varrho(v),x_1)$ follows.

For a detailed proof of part (2), 
 we consider a sequence $x_{1:\infty}$ such
that $x_n$ converges to infinity and
 $x_n/(s_n-x_n)$ converges to $v$. For every positive integer $n$, 
we introduce random variables $(T_n,U_n)$
 distributed as $(X_0,U)$ conditionally on $X_{1:n}=x_{1:n}$.
We first show the convergence in probability of $U_n$, 
then the convergence in distribution of $(T_n,U_n)$.

\begin{lemma}\label{l.pru}
With the notations above,
$U_n$ converges to $v$ in probability.
\end{lemma}
\begin{proof}[Proof of lemma \ref{l.pru}]
Lemma~\ref{l5} yields
$$
\mathbb{P}(T_n=x,U_n\in\mathrm{d} u)=c_np_x\varrho(u)^xb_n(u)q_n(u)\,\mathrm{d} u,
$$
where $c_n$ denotes a normalizing constant which is independent on $x$ and $u$,
$p_x$ depends only on $x$ and $x_1$, $b_n(u)$ depends only on $u$,
$x_n$ and $s_n$, and $q_n(u)$ depends only on $u$ and $n$.
More precisely,
for every integer $x$ such that $x_1\le2x\le2x_1$ and every real number $u$ in $]0,1[$, 
\begin{eqnarray*}
p_x & = & \binom{2x}{x}\binom{x}{x_1-x},
\\
b_n(u) & = & u^{x_n-1/2}(1-u)^{s_n-2x_n-1/2},
\\
q_n(u) & = & 
\sqrt{\frac{(1+u)^n-1}{u}}.
\end{eqnarray*}
We aim to show that, for every integer $x$ such that $p_x$ is positive and
every positive real number $z$,
when $n$ converges to infinity, 
$$
\mathbb{P}(T_n=x,|U_n-v|\ge z)\ll \mathbb{P}(T_n=x).
$$
Since the function $q_n$ is nondecreasing, 
$$
\mathbb{P}(T_n=x,|U_n-v|\ge z)\le c_np_xq_n(1)\int_{|u-v|\ge
  z}\varrho(u)^xb_n(u)\,\mathbf{1}_{[0,1]}(u)\,\mathrm{d} u,
$$
and
$$
\mathbb{P}(T_n=x)\ge c_np_xq_n(0)\int_0^1\varrho(u)^xb_n(u)\,\mathrm{d} u.
$$
The ratio of the two integrals written above
is
$\mathbb{P}(|B_n-v|\ge z)$, where $B_n$ is a beta random variable of
parameters $(\alpha_n,\beta_n)$, with 
$$
\alpha_n=x_n-x+1/2,\quad
\beta_n=s_n-2x_n+2x+1/2.
$$
Since $\alpha_n$ and $\beta_n$ both converge to infinity and
$\alpha_n/(\alpha_n+\beta_n)$ converges to $v$, it is an easy matter
to show that $B_n$ converges in probability to $v$. However, we need a
stronger statement, namely the fact that 
$\mathbb{P}(|B_n-v|\ge z)\ll q_n(0)/q_n(1)$.
Note that $q_n(0)=\sqrt{n}$ and $q_n(1)\sim2^{n/2}$, hence $q_n(0)/q_n(1)\ll1$.

One can write an elementary proof of this, based on the
representation of beta random variables with integer parameters
as ratios of sums of
i.i.d.\ exponential random variables and on large deviations properties of
these sums. Instead, we rely on approximations of beta distributions by
normal distributions provided by Alfers and Dinges~(1984).  A rephrasing of
corollary~1 on page 405 of this paper is as follows.  Let $(Y_k)_k$
denote a sequence of beta random variables of parameters
$(ka_k,k(1-a_k))$. Assume that $k$ converges to infinity and that
$a_k$ converges to a limit $0<a<1$.  Then, for every fixed $y$ such
that $a<y<1$, the ratio
$$
\frac{\mathbb{P}(Y_k\ge y)}{\mathbb{P}\left(Z\ge\sqrt{2k\ell(a_k,y)}\right)}
$$ 
converges to a
finite and positive limit, which depends on $a$ and $y$ only, where $Z$
denotes a standard Gaussian random variable, and $\ell$ denotes the function defined by
$$
\ell(\alpha,y)=\alpha\log\left(\frac{\alpha}{y}\right)
+
(1-\alpha)\log\left(\frac{1-\alpha}{1-y}\right).
$$
Since $a_k$ converges to $a$ and $\ell(\alpha,y)$ is a continuous
function of $\alpha$, standard estimates of Gaussian tails and the result by
Alfers and Dinges
show that
there exists a positive constant $C<1$, independent on $k$, such that
for every $k$ large enough,
$$
\mathbb{P}(Y_k\ge y)\le C^k.
$$
Applying this to our setting, first to the random variables $B_n$ and
to $y=v+z$, then to the random variables $1-B_n$ and to $y=1-v+z$, one
gets the existence of a constant $C<1$ such that, for every $n$ large enough,
$$
\mathbb{P}(|B_n-v|\ge z)\le 2C^{\alpha_n+\beta_n}.
$$
Since $\alpha_n+\beta_n=s_{n-1}+x+1\ge s_{n-1}\gg n$,
$2C^{\alpha_n+\beta_n}\ll q_n(0)/q_n(1)$, and the proof of
lemma~\ref{l.pru} is complete.
\end{proof}

We now apply lemma~\ref{l.pru} to the proof of part~(2).
Introduce the finite sums
$$
p(u)=\sum_xp_x\varrho(u)^x.
$$
For every $u$ in $]0,1[$, the distribution of $T_n$ conditionally on
$U_n=u$ is independent on $n$ and such that
$$
\mathbb{P}(T_n=x\,|\,U_n=u)=p(u)^{-1}p_x\varrho(u)^x.
$$
Hence, for every measurable subset $B$ of $]0,1[$,
$$
\mathbb{P}(T_n=x,U_n\in B)=
\mathbb{E}\left(p(U_n)^{-1}p_x\mathbf{1}_B(U_n)\varrho(U_n)^x\right).
$$  
The function $u\mapsto p(u)^{-1}p_x\mathbf{1}_B(u)\varrho(u)^x$ is bounded
by $1$ on $]0,1[$ and, as soon as $v$ is not in the boundary of $B$,
continuous at $u=v$.  Since $U_n$ converges in distribution to $v$,
this implies that $\mathbb{P}(T_n=x,U_n\in B)$ converges to
$p(v)^{-1}p_x\mathbf{1}_B(v)\varrho(v)^x$, for instance for every interval
$B=[0,u]$ with $u\ne v$.  This is equivalent to the desired
convergence in distribution.

\subsection{Remarks}

For every positive integer $n$ and every admissible
sample, $s_n\ge2x_n(1-1/2^n)$ since $x_k\ge x_{k+1}/2$ for every nonnegative
integer $k$, 
hence $s_n-x_n\ge
x_n+o(x_n)$ and $u\le1$ in the asymptotics that we consider.
Furthermore, the function $\varrho$ decreases from
$\varrho(0^+)=+\infty$ to $\varrho(1^-)=0$.

The measures $\mu(r,x)$ for the first values of $x$ are as follows:
$\mu(r,1)=\delta_1$,
$$
\mu(r,2)=\frac{\delta_1+3r\delta_2}{1+3r},\quad
\mu(r,3)=\frac{3\delta_2+5r\delta_3}{3+5r},\quad
\mu(r,4)=\frac{3\delta_2+30r\delta_3+35r^2\delta_4}{3+30r+35r^2},
$$
and
$$
\mu(r,5)=\frac{15\delta_3+70r\delta_4+63r^2\delta_5}{15+70r+63r^2}.
$$

\subsection{Proof of proposition~\ref{lo}}

The monotonicity with respect to $r$ is valid in a wider setting,
described in proposition~\ref{pw} below, but the monotonicity with
respect to $x$ is more specific.

\begin{proposition}\label{pw}
  Let $\mu$ denote a nonzero bounded measure with exponential moments. For
  every real number $a$, introduce the measures $\nu_a$ and $\mu_a$
  defined by the relations $\nu_a(\mathrm{d}
  x)=\mathrm{e}^{ax}\mu(\mathrm{d} x)$ and $\mu_a=\nu_a/|\nu_a|$. Then
  the family $(\mu_a)_{a}$ is stochastically nondecreasing.
\end{proposition}

\begin{proof}[Proof of proposition \ref{pw}]
  Fix $x$. The derivative of $\mu_a([x,+\infty))$ with respect to $a$
  has the sign of $D(x)$, with
$$
D(x)=\int_{y\ge x}y\mathrm{e}^{ay}\mu(\mathrm{d} y)\,\int\mathrm{e}^{az}\mu(\mathrm{d} z)
-\int_{y\ge x}\mathrm{e}^{ay}\mu(\mathrm{d} y)\,\int z\mathrm{e}^{az}\mu(\mathrm{d} z).
$$
The variations of $D(x)$ with respect to $x$ are
given by
$$
\mathrm{d} D(x)=\mathrm{e}^{ax}\mu(\mathrm{d} x)\,\int(y-x)\mathrm{e}^{ay}\mu(\mathrm{d} y).
$$
The integral in the right hand side is a nonincreasing function of $x$.
Since $D(0)=D(\infty)=0$, the function $x\mapsto D(x)$ is
nondecreasing for $x\le x_a$ and nonincreasing for $x\ge x_a$, 
where $x_a$ solves the equation
$$
 \int y\mathrm{e}^{ay}\mu(\mathrm{d} y)
=
x_a\int\mathrm{e}^{ay}\mu(\mathrm{d} y).
$$
This proves that $D(x)\ge0$ for every $x$, hence
$\mu_a([x,+\infty))\le\mu_b([x,+\infty))$ for every $a\le b$.
This concludes
the proof of proposition~\ref{pw}.
\end{proof}

We turn to the monotonicity of $\mu(r,x)$ with respect to $x$.
We fix a value of $r$ and write every $\nu(r,x)$ as
$$
\nu(r,x)=\sum_ya_y^x\delta_y.
$$
We want to prove that for every $x$, $G(z)\ge0$ for every $z$, with
$$
G(z)=\sum_ya_y^x \sum_{y\ge z}a_y^{x+1}- \sum_{y\ge z}a_y^x \sum_{y}a_y^{x+1}.
$$
One sees that $G(0)=G(\infty)=0$, and simple computations show that
$$
F(z)=\frac1{a_z^{x}}(G(z+1)-G(z))
=\sum_{y}a_y^{x+1}-\frac{a_z^{x+1}}{a_z^x}\sum_ya_y^x.
$$
At this point, we use the specific form of the coefficients $a_z^x$, which yields
$$
\frac{a_z^{x+1}}{a_z^x}=2\frac{(2x+1)(2z-x)}{(x+1)(x+1-z)}.
$$
This shows that $(F(z))_z$ is a nonincreasing
sequence, hence $G(z+1)-G(z)\ge0$ if $z<z_*$
and $G(z+1)-G(z)\le0$ if $z\ge z_*$, for a given $z_*$.  
Hence the sequence $(G(z))_z$ is
nondecreasing on $z\le z_*$ and nonincreasing on $z\ge z_*$. Since
$G(0)=G(\infty)=0$, this implies that 
$G(z)\ge0$ for every positive $z$. This concludes
the proof of proposition \ref{lo}.

%
\section{Limit posterior distributions of initial populations}
\label{s.tb}

\subsection{Expectations}

Let $u$ in $]0,1[$ and $r=\varrho(u)$.
We are interested in the limit
as $x\to\infty$ of the sequence
$$
\frac1{x}\,\mathbb{E}(\xi_x)=\frac{A(r,x)}{x\,B(r,x)},
$$
with the notations
$$
A(r,x)=\sum_yy\,\nu(r,x)(y)=\sum_yy\,\binom{2y}{y}\,\binom{y}{x-y}\,r^y,
$$
and
$$
B(r,x)=|\nu(r,x)|=\sum_y\binom{2y}{y}\,\binom{y}{x-y}\,r^y.
$$

\begin{definition}
For every positive $\lambda$ and $r$,
introduce
$$
C_\lambda(r,z)=(1-4rz(1+z))^{-\lambda}=\sum_{x\ge0}c_\lambda(r,x)\,z^x.
$$
\end{definition}

Starting from the expansion
$$
(1-4z)^{-1/2}=\sum_{x\ge0}\binom{2x}{x}\,z^x,
$$ 
one can write $B(r,x)$ as the coefficient of $z^x$ in the expansion of
$C_{1/2}(r,x)$ along the powers of $z$, namely, $B(r,x)=c_{1/2}(r,x).$
Likewise, $A(r,x)$ is $r$ times the derivative of $B(r,x)$ with
respect to $r$, hence $A(r,x)$ is the coefficient of $z^x$ in the
expansion of $2rz(1+z)C_{3/2}(r,x)$ along the powers of $z$.  This
yields
$$
A(r,x)=2r\,\left(c_{3/2}(r,x-1)+c_{3/2}(r,x-2)\right),
$$
and
$$
x\,B(r,x)=2r\,\left(c_{3/2}(r,x-1)+2c_{3/2}(r,x-2)\right).
$$
\begin{definition}\label{d.gammam}
For every positive $r$, introduce
$$
\gamma(r)=\frac12\left(\sqrt{\frac{1+r}{r}}-1\right),
\qquad
m(r)=\frac{1+\gamma(r)}{1+2\gamma(r)}
=
\frac12\left(1+\sqrt{\frac{r}{1+r}}\right).
$$
\end{definition}
Note that, for every $u$ in $]0,1[$,
$$
\gamma(\varrho(u))=\frac{u}{1-u},\qquad
m(\varrho(u))=\frac1{1+u}=m_u.
$$
\begin{lemma}\label{l.cl}
For every positive $\lambda$ and $r$, when
$x$ converges to infinity,
$$
c_\lambda(r,x)\sim c_\lambda(r)\,x^{\lambda-1}\,\gamma(r)^{-x},
\qquad
c_\lambda(r)=m(r)^\lambda/\Gamma(\lambda).
$$
\end{lemma}

\begin{proof}[Proof of lemma \ref{l.cl}]
This
is a consequence of known expansions of powers of $1/(1-z)$.
First, recall that
$$
(1-z)^{-\lambda}
=
\sum_{x\ge0}d_\lambda(x)\,z^x,\quad
d_\lambda(x)
=
\frac{\Gamma(x+\lambda)}{\Gamma(x+1)\Gamma(\lambda)}
\sim
\frac{x^{\lambda-1}}{\Gamma(\lambda)}.
$$
We use this and the decomposition
$$
1-4rz(1+z)=\left(1-\frac{z}{\gamma(r)}\right)\,\left(1+\frac{z}{\gamma(r)+1}\right),
$$
to get the expansion
$$
C_\lambda(r,z)=
\sum_xd_\lambda(x)\,\left(\frac{z}{\gamma(r)}\right)^x\,
\sum_xd_\lambda(x)\,\left(\frac{-z}{1+\gamma(r)}\right)^x,
$$
which implies
$$
c_\lambda(r,x)=d_\lambda(x)\,\gamma(r)^{-x}\,
\sum_{y=0}^x\left(\frac{-\gamma(r)}{\gamma(r)+1}\right)^y\,d_\lambda(y)\,
\frac{d_\lambda(x-y)}{d_\lambda(x)}.
$$
When $x$ converges to infinity, the ratios $d_\lambda(x-y)/d_\lambda(x)$ converge
to $1$, hence, by dominated convergence,
$$
c_\lambda(r,x)\sim 
d_\lambda(x)\,\gamma(r)^{-x}\,
\sum_{y\ge0}\left(\frac{-\gamma(r)}{\gamma(r)+1}\right)^y\,d_\lambda(y)
=
d_\lambda(x)\,\gamma(r)^{-x}\,\left(1+\frac{\gamma(r)}{1+\gamma(r)}\right)^{-\lambda},
$$
where the equality stems from the definition of 
the coefficients $d_\lambda(\cdot)$.
Plugging the equivalent of $d_\lambda(x)$ into this and 
using the fact that $1+\gamma(r)/(1+\gamma(r))=1/m(r)$, one deduces 
lemma~\ref{l.cl}.
\end{proof}

Lemma \ref{l.cl} for $\lambda=\frac32$ yields that,
when $x$ converges to infinity,
there exists a constant $\alpha$, whose value is irrelevant, such that
$$
A(r,x)\sim2r\alpha\,x^{1/2}\,\gamma(r)^{-x}\,\gamma(r)\,(1+\gamma(r)),$$
and
$$
x\,B(r,x)-A(r,x)\sim2r\alpha\,x^{1/2}\,\gamma(r)^{-x}\,\gamma(r)^2.
$$
Hence $(x\,B(r,x)-A(r,x))/A(r,x)$ converges to $\gamma(r)/(1+\gamma(r))$, and
$$
\frac{A(r,x)}{x\,B(r,x)}
\quad\mbox{converges to}\quad
\frac{1+\gamma(r)}{1+2\gamma(r)}=m(r).
$$
This is the desired convergence of the expectations because, as
mentioned above, the relation $r=\varrho(u)$ means that $m(r)=m_u$.

\subsection{Modes}

To study the mode of $\xi_x$, one compares $\nu(r,x)(y+1)$ to $\nu(r,x)(y)$.
The ratios
$$
\frac{\nu(r,x)(y+1)}{\nu(r,x)(y)}
=
\frac{(y+1/2)\,(x-y)\,r}{(y+1-x/2)\,(y+1/2-x/2)}
$$
 are the terms of a nonincreasing sequence indexed by $y$. Writing
$y$ as $y=x\,(1+s)/(2s)$ with $s\ge 1$, when $x$ is large, one gets
$$
\frac{\nu(r,x)(y+1)}{\nu(r,x)(y)}\sim\,r\,(s^2-1).
$$
This implies that the sequence $(\nu(r,x)(y))_{y}$ is increasing on
$y\le y_*$ and decreasing on $y\ge y_*$, for a value of $y_*$ such
that
$y_*=x\,(1+s_*)/(2s_*)+o(x)$ with $s_*^2=1+1/r.$
Finally, this shows that, when $r=\varrho(u)$, the mode of 
$\mu(r,x)$ is at $x/(1+u)+o(x)$.

\subsection{Distributions}

Our next computation is based on characteristic functions. Fix $u$ in
$]0,1[$ and let $r=\varrho(u)$. For every positive integer $x$, introduce
$$
F_x(t)=\mathbb{E}\left(\exp\left(t\,\frac{\xi_x-x\,m}{\sqrt{x}}\right)\right).
$$
Recall that
$$
m(r)=\frac{1+\gamma(r)}{1+2\gamma(r)},\qquad
\frac1{1+2\gamma(r)}=\sqrt{\frac{r}{1+r}},\qquad 
r=\frac{(1-u)^2}{4u}.
$$
Since
$\mathbb{E}(\exp(t\xi_x))
=
B(r\mathrm{e}^t,x)/B(r,x)$,
$$
F_x(t)=\mathrm{e}^{-t\,\sqrt{x}\,m}\,B(r\mathrm{e}^{t/\sqrt{x}},x)/B(r,x).
$$
We turn to the study of the sequence of functions $(B(\cdot,x))_{x\ge1}$.

Since $B(r,x)=c_{1/2}(r,x)$, a consequence of lemma \ref{l.cl}
is that, when $x$ converges to infinity,
$$
\frac{B(r\mathrm{e}^{t/\sqrt{x}},x)}{B(r,x)}
\sim
\left(\frac{\gamma(r)}{\gamma(r\mathrm{e}^{t/\sqrt{x}})}\right)^{x}
\frac{S(x,t/\sqrt x)}{m(r)^{1/2}},
$$
where, for every $s$,
$$
S(x,s)=
\sum_{y=0}^x\left(\frac{-\gamma(r\mathrm{e}^{s})}
{\gamma(r\mathrm{e}^{s})+1}\right)^y\,d_{1/2}(y)\,
\frac{d_{1/2}(x-y)}{d_{1/2}(x)}.
$$
We get rid of the fraction involving $S(x,t/\sqrt x)$
through lemmas \ref{ld} and \ref{ls}. 
\begin{lemma}\label{ld}
For every nonnegative $x$ and $y$, $d_{1/2}(y)\,d_{1/2}(x)\le d_{1/2}(x+y)$.
\end{lemma}

\begin{proof}[Proof of lemma \ref{ld}]
  A probabilistic proof is as follows.  For every nonnegative $x$,
  $d_{1/2}(x)=2^{-2x}\binom{2x}{x}$ is the probability that a simple
  symmetric random walk on the integer line is at its starting point
  after $2x$ steps.  Hence $d_{1/2}(x+y)$ is the probability that the
  random walk is at its starting point after $2x+2y$ points and
  $d_{1/2}(y)\,d_{1/2}(x)$ is the probability that the random walk is
  at its starting point after $2x$ steps and also after $2x+2y$
  points. The latter event being included in the former, this shows
  the desired inequality.
\end{proof}

\begin{lemma}\label{ls}
When $x$ converges to infinity, $S(x,t/\sqrt{x})$ converges to 
$m(r)^{1/2}$.
\end{lemma}

\begin{proof}[Proof of lemma \ref{ls}]
Since $S(x,0)$ converges to $S(\infty,0)=m(r)^{1/2}$ when $x$ converges to infinity, 
we show that $S(x,t/\sqrt{x})-S(x,0)$ converges to $0$. By lemma \ref{ld}, 
the ratios of coefficients $d_{1/2}$ involved in $S(x,t/\sqrt{x})$ and 
$S(x,0)$ are bounded by $1$.
Adding terms such that $y\ge x+1$, one gets
$|S(x,t/\sqrt{x})-S(x,0)|\le T(t/\sqrt{x})$, where
$$
T(t/\sqrt{x})=\sum_{y=0}^{+\infty}\left|\left(\frac{\gamma(r\mathrm{e}^{t/\sqrt{x}})}
{\gamma(r\mathrm{e}^{t/\sqrt{x}})+1}\right)^y
-
\left(\frac{\gamma(r)}
{\gamma(r)+1}\right)^y\right|.
$$
All the terms in the sum have the same sign, hence 
$$
T(t/\sqrt{x})=
\left|\sum_{y=0}^{+\infty}\left(\frac{\gamma(r\mathrm{e}^{t/\sqrt{x}})}
{\gamma(r\mathrm{e}^{t/\sqrt{x}})+1}\right)^y
-
\left(\frac{\gamma(r)}
{\gamma(r)+1}\right)^y\right|.
$$
One can compute the sum of each geometric series.
This yields
$$
|S(x,t/\sqrt{x})-S(x,0)|\le T(t/\sqrt{x})
=
\left|\gamma(r\mathrm{e}^{t/\sqrt{x}})-\gamma(r)\right|,
$$
which proves the lemma since $\gamma(\cdot)$ is a continuous function.
\end{proof}

Lemma \ref{ls} shows that
$$
\frac{B(r\mathrm{e}^{t/\sqrt{x}},x)}{B(r,x)}\sim
\left(\frac{\gamma(r\mathrm{e}^{t/\sqrt{x}})}{\gamma(r)}\right)^{-x}.
$$
The rest of the proof is standard.
A Taylor expansion of $\gamma(\cdot)$ around $r$ yields
$$
\gamma(r\mathrm{e}^{t/\sqrt{x}})=
\gamma(r)+(\mathrm{e}^{t/\sqrt{x}}-1)\,\gamma'(r)+
(\mathrm{e}^{t/\sqrt{x}}-1)^2\,\gamma''(r)/2+o((\mathrm{e}^{t/\sqrt{x}}-1)^2).
$$
Using the expansion of $\mathrm{e}^{t/\sqrt{x}}$ along powers of $1/\sqrt{x}$ 
and dividing everything by $\gamma(r)$, one gets
$$
\frac{\gamma(r\mathrm{e}^{t/\sqrt{x}})}{\gamma(r)}=
1+\left(r\,\frac{\gamma'(r)}{\gamma(r)}\right)\,\frac{t}{\sqrt{x}}+
\left(r\,\frac{\gamma'(r)}{\gamma(r)}+
r^2\frac{\gamma''(r)}{\gamma(r)}\right)\,\frac{t^2}{2x}
+o\left(\frac1{x}\right).
$$
Note that
$$
r\frac{\gamma'(r)}{\gamma(r)}=-m(r).
$$
Taking logarithms, writing the ratio of functions $\gamma$ as
$$
(\gamma(r\mathrm{e}^{t/\sqrt{x}})/\gamma(r))^{-x}=
\exp(-x\,\log(\gamma(r\mathrm{e}^{t/\sqrt{x}})/\gamma(r))),
$$
and using the expansion $\log(1+z)=z-z^2/2+o(z^2)$ when $z=o(1)$, one gets
that
$F_x(t)$ is equivalent to the exponential of
$$-tm(r)\sqrt{x}-x\,
\left(-m(r)\,t/\sqrt{x}+(m_2(r)-m(r))\,t^2/(2x)-m(r)^2\,t^2/(2x)+o(1/x)\right),
$$
where
$$
m_2(r)=r^2\frac{\gamma''(r)}{\gamma(r)}.
$$
Finally,
$F_x(t)$ converges to $\mathrm{e}^{\sigma^2(r)\,t^2/2}$, with
$$
\sigma^2(r)=m(r)^2+m(r)-m_2(r).
$$
Using the definitions of $m(r)$ and $m_2(r)$ as functions of
$\gamma(r)$ and its derivatives, one gets
$$
\sigma^2(r)=r\,\left(-r\,\frac{\gamma'(r)}{\gamma(r)}\right)'=r\,m'(r).
$$
Using the formula for $m(r)$ at the beginning of this section, one gets finally
$$
\sigma^2(r)=\frac14\sqrt{\frac{r}{(1+r)^3}}=\frac{u\,(1-u)}{(1+u)^3}
=\sigma^2_u.
$$
The proof is complete.

%
\section{Conditional hitting times}
\label{s.tc}
\subsection{Proof of theorem \ref{tc}}
We introduce the
renewal process $(\zeta_x)_{x\ge1}$ with increments $(\varepsilon_x)_{x\ge1}$, that is
$$
\zeta_x=\inf\{y\ge1\,;\,\sigma_y\ge x\}.
$$
The usual central limit theorem for renewal processes states that
$(\zeta_x-mx)/\sqrt{x}$ converges in distribution to a centered
Gaussian variable whose variance is the variance $u(1-u)$ of every $\varepsilon_x$
divided by the cube of the mean $1+u$ of every $\varepsilon_x$, 
that is $u(1-u)/(1+u)^3=\sigma^2_u$.  

Our next lemma expresses the distribution of $\eta_x$ for every positive $x$
in terms of the 
distributions of the random variables $(\zeta_z)_{1\le z\le x+1}$.

\begin{lemma}\label{l.eta}
For every positive $x$ and $y$,
$$
\mathbb{P}(\eta_x=y)
=
\frac{1+u}{1-(-u)^{x+1}}\,\sum_{z=0}^{x-1}(-u)^z\mathbb{P}(\zeta_{x+1-z}=y+1).
$$
\end{lemma}

\begin{proof}[Proof of lemma \ref{l.eta}]
Let $x$ and $y$ denote positive integers.
We begin with the fact that
$$
\{\zeta_{x+1}=y+1\}=\{\sigma_{y}=x\}\cup\{\sigma_{y}=x-1,\varepsilon_{y+1}=2\},
$$
hence
$$
\mathbb{P}(\sigma_{y}=x)=\mathbb{P}(\zeta_{x+1}=y+1)-u\,\mathbb{P}(\sigma_{y}=x-1).
$$
Iterating this recursion, one gets
$$
\mathbb{P}(H_x)\,\mathbb{P}(\eta_x=y)=\mathbb{P}(\sigma_{y}=x)
=
\sum_{z=0}^{x-1}(-u)^z\,\mathbb{P}(\zeta_{x+1-z}=y+1).
$$
Summing over every positive value of $y$ and using the facts that $\mathbb{P}(\zeta_{z}=1)=0$ if 
$z\ge3$ and that $\mathbb{P}(\zeta_{2}=1)=u$, one gets
$$
\mathbb{P}(H_x)=(-u)^{x-1}(1-u)+\sum_{z=0}^{x-2}(-u)^z=\frac{1-(-u)^{x+1}}{1+u}.
$$
This concludes the proof.
\end{proof}

Lemma \ref{l.eta}, the fact that $|u|<1$
and the convergence of the distribution of
$(\zeta_x-mx)/\sqrt{x}$, imply the same convergence for the distribution of
$(\eta_x-mx)/\sqrt{x}$.

Finally, $(\xi_x-mx)/\sqrt{x}$, $(\zeta_x-mx)/\sqrt{x}$ and
$(\eta_x-mx)/\sqrt{x}$ all converge in distribution to the same limit,
which is the centered Gaussian distribution with variance $\sigma^2_u$.

%
\subsection{Sharp bounds}

\begin{lemma}\label{l.ex}
For every positive $x$,
$$
\mathbb{E}(\eta_x)=\frac{x+1}{1+u}\,
\frac{1+(-u)^{x+2}}{1-(-u)^{x+1}}-\frac{1+u^2}{(1+u)^2}.
$$
\end{lemma}

For instance,
$$
\mathbb{E}(\eta_1)=1,\quad\mathbb{E}(\eta_2)=2-\frac{u}{1-u(1-u)}.
$$
For every positive integer $x$, one can deduce from the exact formula above that
$$
\frac{x}{1+u}-\frac{2u^2}{(1+u)^2}\le
\mathbb{E}(\eta_x)\le\frac{x}{1+u}+\frac{2u}{(1+u)^2}.
$$
The width of the interval delimited by the upper and the lower bounds
of $\mathbb{E}(\eta_x)$ above is $2u/(1+u)\le1$.

Bounds on $\mathbb{E}(\eta_x)$, depending on the parity of $x$,  
are as follows. For every odd $x$,
$$
\mathbb{E}(\eta_x)\ge x/(1+u),
$$
and for every even $x$,
$$
\mathbb{E}(\eta_x)\le x/(1+u)+u(1-u)/(1+u)\le(x+1/4)/(1+u).
$$
These refined bounds yield intervals around $\mathbb{E}(\eta_x)$, which depend on the
parity of $x$, and whose width is always at most $2u/(1+u)^2\le1/2$.

\begin{proof}[Proof of lemma \ref{l.ex}]
Fix a positive integer $x$, a real number $u$ in $]0,1[$, and let $r=\varrho(u)$.
Let $p_x^y=\mathbb{P}(\sigma_y=x)$. Then $\eta_x(\mathbb{P})$ is proportional
to the measure $\displaystyle\sum_yp_x^y\,\delta_y$ and, for every positive $y$,
$$
\sum_{x=y}^{2y}p_x^y\,t^x=\left((1-u)\,t+u\,t^2\right)^y.
$$
Hence the distribution of $\eta_x$ is $\mu_\eta(r,x)$, where
 $\mu_\eta(r,x)=\nu_\eta(r,x)/|\nu_\eta(r,x)|$ 
and
$$
\nu_\eta(r,x)=\sum_{y=h(x)}^x\binom{y}{x-y}\,4^yr^y\,\delta_y.
$$
When $u=0$, $r=\infty$ and $\mu_\eta(\infty,x)$ is the Dirac 
distribution at $x$.
When $u=1$, $r=0$ and $\mu_\eta(0,x)$ is the Dirac distribution at $h(x)$.
For the first values of $x$, the distributions $\mu_\eta(r,x)$ are as follows:
$\mu_\eta(r,1)=\delta_1$,
$$
\mu_\eta(r,2)=\frac{\delta_1+4r\delta_2}{1+4r},\quad
\mu_\eta(r,3)=\frac{2\delta_2+4r\delta_3}{2+4r},\quad
\mu_\eta(r,4)=\frac{\delta_2+12r\delta_3+16r^2\delta_4}{1+12r+16r^2}.
$$
This implies that, for every positive $x$,
$$
\mathbb{E}(\eta_x)=r\,g'_x(r)/g_x(r),\quad 
g_x(r)=\sum_{y=h(x)}^x\binom{y}{x-y}\,4^yr^y.
$$
To study the generating functions $g_x$, we introduce $g_0(r)=1$ and
$$
G(r,z)=\sum_{x\ge0}g_x(r)\,z^x.
$$
Summing first over $y\le x\le 2y$, then over $y\ge0$,
one gets
$$
G(r,z)=\sum_{y\ge0}(4rz)^y\,(1+z)^y=1/(1-4rz(1+z))=C_1(r,z).
$$
From the proof of lemma \ref{l.cl}, one knows that the poles of $C_1(r,z)$ are $z=\gamma(r)$ and $z=-\gamma_2(r)$ with $\gamma_2(r)=\gamma(r)+1$,
hence,
$$
G(r,z)=\frac1{\gamma(r)+\gamma_2(r)}\left(\frac{\gamma_2(r)}{1-z/\gamma(r)}+\frac{\gamma(r)}{1+z/\gamma_2(r)}\right).
$$
This shows that, for every nonnegative $x$,
$$
g_x(r)=\frac{\gamma(r)\gamma_2(r)}{\gamma(r)+\gamma_2(r)}
\left(\gamma(r)^{-(x+1)}-(-\gamma_2(r))^{-(x+1)}\right).
$$
From here, the expression of $\gamma(r)$ as a function of $r$ and
tedious computations of derivatives yield the result.
\end{proof}

%


%

\subsection{Comparison with a na\"\i ve estimator}

For a given value $u$ in $]0,1[$ and for a branching process 
$X_{0:\infty}$ with  offspring distribution
$(1-u)\delta_1+u\delta_2$, when $n$ converges to infinity,
$$
S_n\sim X_n(1+1/(1+u)+1/(1+u)^2+\cdots)=X_n(1+1/u)\quad\mbox{almost surely},
$$ 
hence $B(X_{0:\infty})=u$
almost surely.  The na\"\i ve pointwise
prediction of the mean initial population conditional on $X_1=x$, namely
$N_u(x)=x/(1+u)$, should be compared to
the Bayesian prediction $\mathbb{E}_u(\xi_x)$ for $r=\varrho(u)$.
For $x=2$, one gets
$$
\frac{\mathbb{E}_u(\xi_2)}{N_u(2)}=\frac{(4u+6(1-u)^2)\,(1+u)}{2(4u+3(1-u)^2)}.
$$
This ratio is $1$ when $u=0$ or $u=1$, greater than $1$ for every $u$
in $]0,\frac13[$, and smaller than $1$ for every $u$ in $]\frac13,1[$. Hence
the na\"\i ve and Bayesian predictions cannot be easily compared, at
least on $X_1=x$ for a given finite $x$.

%


\begin{thebibliography}{99}
\footnotesize

\bibitem{alf}
{\sc Alfers, Dieter and Dinges, Hermann} (1984).
A normal 
approximation for beta and gamma tail probabilities.
{\em Zeitschrift f\"ur Wahrscheinlichkeitstheorie und Verwandte Gebiete} 
{\bf 65} (3), 399-420. 

\bibitem{a}
{\sc Jagers, Peter and Klebaner, Fima} (2003). 
Random variation and concentration effects in PCR. 
{\em Journal of Theoretical Biology} {\bf 224}, 299-304.
\\
Available as a preprint at 
\texttt{www.math.chalmers.se/Math/Research/Preprints/2002/98.ps.gz}

\bibitem{kw}
{\sc Kass, Robert E. and Wasserman, Larry A.} (1996). 
The selection of prior
distributions by formal rules. 
{\em Journal of the American Statistical
Association} {\bf 91}, 1343-1370.
\\
Available at 
\texttt{www.stat.cmu.edu/\%7ekass/papers/rules.pdf}


\bibitem{ks}
{\sc Klebaner, Fima C. and Sagitov, Serik} (2002). 
The age of a Galton-Watson population with a geometric 
offspring distribution.  
{\em Journal of Applied Probability} {\bf 39}, 816--828.


\bibitem{b}
{\sc Lalam, Nadia} (2007).
Statistical inference for quantitative polymerase chain reaction using a hidden Markov model: a Bayesian approach.
{\em Statistical Applications in Genetics and Molecular Biology} {\bf 6} (1), 
Article 10.
\\
Available at
\texttt{www.bepress.com/sagmb/vol6/iss1/art10}

\bibitem{c}
{\sc Lalam, Nadia and Jacob, Christine} (2007).
Bayesian estimation for quantification by real-time polymerase 
chain reaction under a branching process model of the DNA 
molecules amplification process.
{\em Mathematical Population Studies} {\bf 14} (2), 111-129.

\bibitem{d}
{\sc Mendoza, Manuel and Guti\'errez-Pe{\~n}a, Eduardo} (2000).
Bayesian conjugate analysis of the Galton-Watson process.
{\em Test} {\bf 9} (1), 149-171.
\\
Available as a preprint at 
\texttt{allman.rhon.itam.mx/\%7emendoza/Final.ps}

\bibitem{e}
{\sc Molina,  Manuel, Gonz\'alez, Miguel and Mota, Manuel} (1998).
Bayesian inference for bisexual Galton-Watson processes.
{\em Communications in statistics. Theory and methods} {\bf 27} (5), 1055-1070.

\bibitem{f}
{\sc Peccoud, Jean and Jacob, Christine}  (1996).
Theoretical uncertainty of measurements using quantitative 
polymerase chain reaction.
{\em Biophysical Journal} {\bf 71} (1), 101-8.
\\
Available at 
\texttt{www.biophysj.org/cgi/content/abstract/71/1/101}

\bibitem{g}
{\sc Piau, Didier} (2005). 
Confidence intervals for non homogeneous branching processes 
and PCR reactions. {\em The Annals of Probability} {\bf 33}, 674-702.
\\
Available at 
\texttt{arxiv.org/abs/math/0503659}

\bibitem{h}
{\sc Piau, Didier} (2002). 
Mutation-replication statistics of polymerase chain reactions. 
{\em Journal of Computational Biology} {\bf 9}, 831-847.
\\
Available as a preprint at 
\texttt{citeseer.ist.psu.edu/491286.html}

\bibitem{i}
{\sc Prakasa Rao, B.L.S.} (1992).
Nonparametric estimation for Galton-Watson type
process. 
{\em Statistics and Probability Letters} {\bf 13}, 287-293.

\bibitem{j}
{\sc Scott, David} (1987).
On posterior asymptotic normality and asymptotic normality of
estimators for the Galton-Watson process. {\em Journal of the
Royal Statistical Society, Series B (Methodological)} {\bf 49} (2), 209-214.

\bibitem{k}
{\sc Sun, Fengzhu} (1995).
The polymerase chain reaction and branching processes. 
{\em Journal of Computational Biology} {\bf 23}, 3034-3040.
\\
Available 
at \texttt{www-rcf.usc.edu/~fsun/Publication/PCR/pcr.pdf}

\bibitem{l}
{\sc Weiss, Gunter  and  von Haeseler, Arndt} (1997).
A coalescent approach to the polymerase chain reaction.
{\em Nucleic Acids Research} {\bf 25} (15), 3082-7.
\\
Available at
\texttt{nar.oxfordjournals.org/cgi/content/abstract/25/15/3082}

\bibitem{m}
{\sc Weiss, Gunter  and  von Haeseler, Arndt} (1995).
Modeling the polymerase chain reaction.
{\em Journal of Computational Biology} {\bf 2} (1), 49-61.

\end{thebibliography}
\end{document}